\documentstyle[12pt]{article}
\newfont{\mib}{cmmib10 scaled 1200}

\renewcommand{\phi}{\varphi}
\newtheorem{Lemma}{Lemma}[section]

\newtheorem{Proposition}{Proposition}[section]

\newtheorem{Definition}{Definition}[section]

\newtheorem{Remark}{Remark}[section]
\newtheorem{Conjecture}{Conjecture}[section]
\begin{document}
\title{Boundedness theorem for Fano log-threefolds}
\author{Alexandr Borisov\\
Department of Mathematics \\
Pennsylvania State University\\
e-mail: borisov@math.psu.edu}
\date{(Revised English version)\\Jan 6, 1994}
\maketitle

\begin{abstract}

 The main purpose of this article is to prove that the family of all Fano
threefolds with log-terminal singularities with bounded index is bounded.
\end{abstract}
\section{Introduction}
First of all, let me recall necessary definitions and list some known results
and conjectures in direction of bondedness of Fano manifolds. All varieties in
this paper are over field of complex numbers.

\begin{Definition}
Normal variety $X$ is called  the variety with log-terminal singularities
if  $mK_X$ is a Cartier divisor for some integer $m$ and there exists a
resolution $\pi:Y\longrightarrow X$ of singularities of $X$ such that
exceptional divisors $F_i$ of $\pi$ have simple normal crossings and in
formula  $K_Y=\pi^*K_X+\sum(a_iF_i)$ all  $a_i>-1$
\end{Definition}

\begin{Definition}
Index (or Gorenstein index) of variety X is a minimal natural number m, s.t.
$mK_X$ is a Cartier divisor. Of course, index is defined for $Q$-Gorenstein
varieties only.
\end{Definition}

\begin{Definition}
Three-dimensional algebraic variety $X$ is called Fano log-threefold if the
following conditions hold.

$1)$ $X$ has log-terminal singularities,

$2)$ $X$ is Q-factorial,

$3)$ Picard number $\rho(X)=1$,

$4)$ $-K_X$ is ample.

\end{Definition}
\begin{Remark} This is only my terminology inspired by the term "Q-Fano
threefolds".
\end{Remark}
The following statement will be proven  in this paper.

{\bf Main Theorem.

 For an arbitrary natural $n$ all Fano log-threefolds of index $n$ lie in
finite number of families.}

\begin{Remark} Unfortunately, no effective bound on any invariant of $X$ will
be given because of Noetherian induction in the section $4.$
\end{Remark}
Here are some results in the direction of boundedness of Fano manifolds.

1) Boundedness theorem for smooth Fano manifolds of an arbitrary direction
is proven  by Koll\'ar, Miyaoka and Mori in \cite{KoMiMo1}. Before this result
there were several proofs with  extra condition $\rho(X)=1$. Three-dimensional
smooth case was also treated before by a long work of many authors beginning
with Fano itself. See \cite{Isk} for discussion.

2)Two-dimensional Fano varieties are traditionally called Del Pezzo surfaces.
Smooth (=terminal) case is fairly easy and the answer is the following.

$P^1\times{P^1}, P^2$ with $0$ to $8$ blown-up points in general position. (The
generality of position may be stated precisely.)

I should notice here that there are many difficult problems concerning Del
Pezzo surfaces if basic field is NOT algebraically closed.

Log-terminal case (with arbitrary Picard number) was studied by Alexeev and
Nikulin (see\cite{Ni}). One of the main results in this direction is
boundedness under the condition of bounded multiplicity of singularities. Let
me mention that by using methods of this paper one can obtain a new simple
proof of some intermediate result, namely boundedness under the condition of
bounded index.

3)The model case of toric varieties of arbitrary dimension is treated in
\cite{BB},
see also \cite{B}.

4) Boundedness of Fano threefolds with $\rho=1$ and terminal singularities
is proven by Kawamata in \cite{Ka1}.

5) Boundedness of Fano threefolds with terminal singularities with no extra
conditions is announced by Mori.

All these results justify the following conjecture.

\begin{Conjecture} The family of all Fano varieties of  given dimension with
discrepancies of singularities greater (or greater or equal) $-1+\epsilon$,
where $\epsilon$ is an arbitrary given positive real number, is bounded.
\end{Conjecture}
\begin{Remark} This conjecture is so natural that probably many people
suspected it but I didn't see it published. Batyrev proposed the weaker variant
of this conjecture, where the condition on discrepancies is replaced by the
condition of boundedness of index. (\cite{Ba2}) Very recently Alexeev told me
that he also stated the above conjecture as a part of a general phenomenon
noticed by Shokurov that some geometric invariants (in this case minimal
discrepancies of Fano varieties) can accumulate only from above (below). See
\cite{Alex} for a discussion.
\end{Remark}
I am expressing my thanks to V. Iskovskih who encouraged me to work in this
direction. I am glad to thank V. Shokurov and V. Alexeev who invited me to the
geometry seminar at Johns Hopkins University and whose remarks simplified and
even corrected this paper. I also want to thank my brother Lev for helpful
discussions.

\section{Preliminary remarks and first lemmas}

In \cite{Ko} Koll\'ar proved that all three-dimensional normal varieties $X$
with an ample Cartier divisor $D$ lie in finite number of families if two
higher coefficients of Hilbert polynomial $P(m)=\chi(mD)$ are bounded. In our
case of three-dimensional Fano varieties of index $n$ it works as follows. Let
$D$
be equal to $-nK_X$. Then it is a Cartier divisor and it follows from general
theory of Riemann-Roch that

 $\chi(O_X(-mnK_X))=\frac1{12}(-K_X)^3nm(nm+1)(2nm+1)+\alpha m+\beta$,
where $\alpha$ and $\beta$ are some constants depending on $X$.

Therefore in order to prove the Main Theorem we only need to prove that
$(-K_X)^3$ is bounded. The following lemma shows that in our case it is also
equivalent to the condition that $h^0(-2nK_X)$ is bounded.

\begin{Lemma}\label{L1}
For arbitrary Fano log-threefold $X$ of index $n$ (actually, only conditions
(1) and (3) are used) the following inequality holds.

$h^0(-2nK_X)\ge{(-K_X)^3 (\frac53 n^3+\frac12 n^2)-1}$
\end{Lemma}

{\bf Proof} By the Kawamata-Vieweg vanishing theorem $h^i(-mnK_X)=0$, $ i>0,
m\ge {0}$. Therefore $h^o(-mnK_X)=\frac1{12}(-K_X)^3nm(nm+1)(2nm+1)+\alpha
m+\beta$ for $m\ge{0}$. Let us consider "the second derivative at 1".

$h^o(-2nK_X)-2h^o(-nK_X)+h^0(O_X)=(-K_X)^3(\frac53 n^3+\frac12 n^2).$

Now the statement of lemma follows from the fact that $h^0(O_X)=1$ and
$h^0(-nK_X)\ge{0}$.

\begin{Lemma}\label{L2}
Suppose $v\in V$ - is a closed point of $k$-dimensional variety with
multiplicity of local ring $r$, $D$ is a semiample $Q$-Cartier divisor on $V$.
Suppose further that the general point $x$ of $V$ can be connected to $v$ by
some curve $\gamma _x$, such that $\gamma _x\cdot D\le d$

Then $D^k\le r\cdot d^k$.
\end{Lemma}

{\bf Proof} For sufficiently large $m$ such that $mD$ is a Cartier divisor one
have that $h^o(O_V(mD))=\frac{m^kD^k}{k!}+O(m^{k-1})$. Therefore if $D^k>r\cdot
d^k$ then for $m>>0$ one can find a non-zero global section $s\in H^0(O_V(mD))$
such that its image by trivialization map of $O_V(mD)$ in $v$ lies in
$(md+1)$th power of maximal ideal of point $v$. Then every curve $\gamma _x$
lies in $Supp(s)$, that is impossible.
\begin{Remark}
The above lemma is very general. In applications $V$ will be our Fano
log-threefold $X$ and $D$ will be $(-K_X)$.
\end{Remark}

\section{Covering family and first division into cases}
\begin{Remark} (about notations). We will often consider birational varieties.
Doing this we will usually identify curves on different varieties if they
coincide in their general points. Namely, let $X\leftarrow - \rightarrow
X^\prime$ and $L\subset X,  L^\prime\subset X^\prime$  be curves. Then $L$ and
$L^\prime$ are identified if there are Zariski open subsets $U\subset X$ and
$U^\prime\subset X^\prime$, such that the above rational map is defined on them
and $U\cong U^\prime$, $L\cap U \cong L^\prime \cap U^\prime\ne\emptyset$ via
it. The identified curves will be usually denoted by the same symbol. The same
convention will be used for two-dimensional subvarieties. If it is necessary to
point out that, say, simple divisor $S$ is considered on variety $X$ it will be
denoted by $S_X$. Another convention is that $\{l\}$ will denote the family of
curves with general element $l$ and $\{H\}$ will denote the LINEAR system of
Weil divisors with general element $H$. It will be clear in every particular
case why these conventions agree with each other.
\end{Remark}

Now we start to prove our Main Theorem. Suppose $X$ is a Fano log-threefold,
$\pi^Y_X :Y\longrightarrow X$ is its $Q$-factorial terminal modification,
$\pi^{Y_1}_Y\longrightarrow Y$ is a resolution of isolated singularities of
$Y$. By the Miyaoka-Mori theorem (\cite{MM}, see also \cite{Ka2}) there exists
a  covering family of rational curves $\{l\}$, such that $l\cdot (-K_X)\le 6$.
The family $\{l\}$ is free on $Y_1$ that is small full deformation of $l$
covers small neighborhood of $l$. ( See \cite{Na}.) We can and will denote by
$\{l_{Y_1}\}$ full family, that is (some Zariski open subset of) a component of
the scheme of morphisms from $P^1$ to $Y_1$. Consider the RC-fibration $\phi :
Y_1-{ - }\rightarrow Z$, associated with $\{l\}$. (See \cite{KoMiMo1},
\cite{KoMiMo2}.) The following cases are possible.

$(0)$ $dimZ=0$. In \cite{KoMiMo2} such $X$ are called primitive. It implies
that two general points of $Y_1$ can be joined by chain of no more than $3$
curves from $\{l\}$. It follows now from one of the "gluing lemmas"
(\cite{KoMiMo2}) that we can glue them together and obtain new family
$\{l^{\prime}\}$. Then we can apply to it lemma \ref{L2} and obtain that

$(-K_X)^3 \le (3\cdot l\cdot (-K_X))^3 \le (3\cdot 6)^3$

$(1)$ $dimZ=1$. In this case after some additional blowing-up $\tilde{Y}
\longrightarrow Y_1$ we obtain a morphism $\phi _{\tilde{Y}}:\tilde{Y}
\longrightarrow Z$. Here $Z\cong{P^1}$, because $X$ is rationally connected
(see \cite{KoMiMo2}).

$(2)$ $dimZ=2$. In this case general $l\in \{l\}$ is smooth and does not
intersect with another general $l$ on $Y_1$. And it is exactly the general
fiber of the RC-fibration.

We will proceed  by the following way. First of all we will treat the case (1).
Doing this we will require $l\cdot (-K_X)$ to be bounded not by $6$ but only by
an arbitrary constant depending on $n$. After that we will reduce the case (2)
to the case (1) but for some new family $\{l^{\prime}\}$ where $l^{\prime}\cdot
(-K_X)$ will be bounded.

\section{The treatment of case (1)}
Let $S$ be a general fiber of our RC-fibration. As we already mentioned, the
image of RC-fibration is rational. This implies that S are linear equivalent on
$Y_2$ and therefore on $X$. Notice that it can happen that $\{l\}$ does not
connect two general points of $S$ immediately. But it will always be true if we
glue two examples of $\{l\}$. (See \cite{KoMiMo2}.) Therefore we will assume,
that $\{l\}$ is a connecting family on $S$. Evidently, $l^2\ge1$ on a smooth
surface $\tilde{S}=S_{\tilde{Y}}$. The condition that $X$ is Q-factorial with
Picard number $1$ implies that $S_X=\alpha H, \alpha >0$, where $H=(-2nK_X)$.
We will assume up to the end of this section that $l\cdot H\le \rho,$ where
$\rho$ is some constant depending on $n.$
\begin{Proposition}\label{S1}
If $h^0(H)>2(\rho +1)^2$ then $\alpha\le\frac12$.
\end{Proposition}
{\bf Proof} Let $S_1$ and $S_2$ be two general surfaces from $\{S_X\}$.
Let $l_1, l_2,...l_{\rho+1} \subset S_1$ and $l_{\rho+2}, l_{\rho+3},...
l_{2\rho+2} \subset S_2$ be general curves from $\{l\}$. We have that $H\cdot
l\le \rho$, therefore
$$dimH^0(O_X(H))-dimH^0(J_{l_i}\cdot O_X(H))\le dimH^0(O_{l_i}(H))\le\rho +1.$$
(Here $J_{l_i} $ is an ideal sheaf of the curve $l_i\subset X$.)

 This implies that $$codim(\bigcap^{2\rho +2}_{i=1}H^0(J_{l_i}\cdot
O_X(H)))\le\sum^{2\rho +2}_{i=1}(\rho+1)=2(\rho+1)^2.$$
If $h^0(H)>2(\rho+1)^2$ then there exists a divisor $H^*\in |H|$, such that all
$l_i\subset H^*$. Suppose $\pi_1$ is a composition of birational morphism
$\tilde{S_1}\longrightarrow S_1$ and embedding $S_1\longrightarrow X$. If $H^*$
does not contain $S_1$ then $(\pi^{\tilde{Y}}_Y)^* (H^*)$ does not contain
$\tilde{S}_1$ but at the same time contains preimages of $l_i$. On
$\tilde{S}_1$ we have $(\pi^*_1H^*)\cdot l\ge\sum l_i\cdot l\ge\rho+1$. It
contradicts to the fact that $(\pi_1^*(H^*)\cdot l=H\cdot l\le\rho$. Therefore
$H^*$ contains $S_1$ and, by the same arguments, $S_2$. This implies that
$\alpha \le \frac12.$

We will always assume below that $\alpha \le \frac12.$
\begin{Proposition}\label{S2}
For arbitrary $S\in \{S\}$ on $X$, arbitrary positive integer $k$

$h^i(X,J_S\cdot O_X(kH))=0$ for every $i>O$.
\end{Proposition}
{\bf Proof} $S$ is a simple divisor, therefore $J_S=O_X(-S)$, where $O_X(-S)$
is a divisorial sheaf sheaf, associated with Weil divisor $(-S)$. After that,
$O_X(kH)$ is an invertible sheaf, therefore $J_S\cdot O_X(kH) =O_X(kH-S)$. Now
one can apply Kawamata-Vieweg vanishing theorem (see the reformulation of it in
\cite{Al}) because $kH-S-K_X=(k-\alpha+\frac1{2n})H$ is ample for $k \ge 1.$
\begin{Proposition}\label{S3}
For all $k>0$, $l>0$  $h^i(S,O_S(kH))=O.$
\end{Proposition}
{\bf Proof} It follows from exact sequence

$0 \longrightarrow J_S \cdot O_X(kH) \longrightarrow O_X(kH) \longrightarrow
O_S(kH) \longrightarrow 0$,
vanishing theorem and proposition \ref{S2}.
\begin{Proposition}\label{S4}
All surfaces $S_X$ for given $n$ and $\rho = c(n)$ lie in finite number of
families.
\end{Proposition}
{\bf Proof} By a result of Koll\'ar (\cite{Ko}) and proposition \ref{S3} it is
enough to prove the boundedness of coefficients of Hilbert polynomial
$P(k)=\chi(O_S(kH))=h^0(O_S(kH)), k\ge 1.$ For this purpose we will prove that
there exists some constant $c_1(n,\rho )$ such that for all $k\ge 1$ it is true
that $h^0(O_S(kH))\le k^2\cdot c_1(n,\rho )$. It implies the boundedness of
coefficients by the following arguments. Suppose $P(k)=a_2k^2+a_1k+a_0$.
Evidently, $0\le a_2\le c_1$. Therefore $|a_1|=|P(2)-P(1)-3a_2|\le 4c_1$. After
that, $|a_0|=|P(1)-a_2-a_1|\le 5c_1$.

In order to prove that $h^0(O_S(kH))\le k^2\cdot c_1(n,\rho )$ consider the
following construction. By applying several times gluing lemma to a free family
$\{l\}$ on $\tilde{S}$ (\cite{KoMiMo2}) we obtain families $\{l_k\}$ such that
$l_k=k\cdot l$ as divisors on $\tilde{S}$ and therefore on $S$. (Here "=" means
algebraic equivalence.) Notice that the natural map $\mu _k
:H^0(S,O_S(kH))\longrightarrow H^0(l_{k\rho +1},O_{l_{k\rho +1}}(kH))$ is
injective. Otherwise there should have been some $D\in |kH|$ containing
$l_{k\rho +1}$ but not containing $S$. As in the proof of proposition \ref{S1}
we obtain a contradiction by intersecting with general $l$. The fact that
$\mu_k$ is injective implies that $h^0(O_S(kH))\le h^0(l_{k\rho +1},O_{l_{k\rho
+1}}(kH))\le l_{k\rho +1}(kH)+1=(k\rho +1)k\rho +1$, that is what we need.
\begin{Proposition}\label{S5}
In the condition of the above proposition there is a constant $c_2(n,\rho )$,
such that on every general $S_X$ EVERY two points can be joined by some
irreducible curve $\gamma$, such that $\gamma\cdot (-K_X)\le c_2$.
\end{Proposition}
{\bf Proof} It is a straightforward consequence of boundedness of $S_X$ with
$H|_{S_X}$. Indeed, it is true for a general element of every one of families
in proposition \ref{S4} and Noetherian induction on base completes the proof.

\begin{Remark} Of course, two GENERAL points of $S_X$ are already connected by
$l$, but the above proposition gives much more.
\end{Remark}
Now we can complete the treatment of the case $(1)$. By the definition of Fano
log-threefold $\rho (X)=1$ therefore two general $S_X$ intersect with each
other. Moreover, they intersect along some curve $C$ because $X$ is
Q-factorial. We know that $\{S_X\}$ is a linear system, therefore all of them
contain $C$. It may happen that $C$ lies in $Sing(X)$, but the multiplicity of
$X$ in a general point $x_0\in C$ is bounded by $2n$, because the index of $X$
is bounded by $n$. (By canonical cover trick it is a factor of $CDV$
singularity that is analytically isomorphic to $(DV-point)\times (disk).$)
Therefore we can apply lemma \ref{L2} to $X, (-K_X), x_0$  to obtain a bound on
$(-K_X)^3.$

\section{Two lemmas}
In this section we will prove some adjunction lemma and a lemma about accurate
resolution that will be used in next section to treat the case $(2)$. However,
these lemmas themself are interesting enough to deserve  a separate section.
\begin{Lemma}\label{LA}(adjunction)
Suppose $X$ is a three-dimensional variety and $S$ is simple Weil divisor on
it, such that $(K_X+S)$ is Q-Cartier. Suppose $\{L\}$ is a covering family of
curves on $S$, $\hat{S}$ is a minimal resolution of normalization of $S$. Then
$K_{\hat{S}}\cdot L\le (K_X+S)\cdot L.$
\end{Lemma}
{\bf Proof} Denote by $\pi$ the natural morphism $\hat{S}\longrightarrow X.$
Then by the proposition 3.2.2 of \cite{SH} $K_{\hat{S}}=\pi ^*(K_X+S)-D$, where
D is an effective divisor. The rest is trivial.
\begin{Remark}
The above lemma is due to Shokurov. In the first variant of this paper I
formulated and proved it only under condition that singularities of $X$ were
isolated which is enough for applications.
\end{Remark}
\begin{Lemma}\label{LR}(accurate resolution)

Suppose $X$ is a Q-factorial three-dimensional variety, $E\subset X$ is a
simple Weil divisor, $\{L\}$ is a covering family of curves on E. Suppose
further that there exists a covering family $\{l\}$ on $X$, such that $l\cdot
E\ge 1$ and a linear system $|H|$ on $X$, such that the following inequalities
hold true. ($c_i$ are some nonnegative constants.)

$1)$ $H\cdot l\le c_1$

$2)$ $H\cdot L\le c_2$

$3)$ $K_X\cdot L\le c_3$

$4)$ $-E\cdot L\le c_4$

Then $h^0(H)>1+(c_1+1)(c_2+c_1c_4+1)$ implies that there exists a resolution
$Y\longrightarrow X$, such that $\{L\}$ have no base points on $E_Y$ and

$K_Y\cdot L\le c_3+2(c_2+c_1c_4).$
\end{Lemma}
\begin{Remark}
The proof of this lemma will be pretty long. It will take the rest of the
section.
\end{Remark}
\begin{Remark}
In some sense this lemma is a very weak substitute for the following conjecture
for which I have a lot of evidence.
\end{Remark}
{\bf Accurate Resolution Conjecture} For an arbitrary Q-Gorenstein threefold
$X$ there exists a resolution of singularities $\pi:Y\longrightarrow X,$ such
that for EVERY Q-Cartier divisor $H$ on $X$ containing a curve $L_X$ not lying
in Sing(X) the following inequality holds true.

$(K_Y+D_Y)\cdot L_Y \le (K_X+D_X)\cdot L_X$

 $ $

First of all we will introduce some convenient notations. Let $\{D\}$ be a
linear system of Weil divisors. We will denote by $H^0(\{D\})$ the
corresponding vector subspace in $H^0(O_X(D))$, where $O_X(D)$ is a divisorial
sheaf, associated with $D$. Reversely, for a linear subspace $V\subset
H^0(\{D\})$ let $|V|$ be the corresponding linear system. Divisor that
corresponds to $s\in H^0(O_X(D))$ will be denoted by $(s)$. Section that
determines divisor $D$ will be called "equation" of $D.$ Of course, it is
defined up to multiplicative constant. By definition
$h^0(\{D\})=dimH^0(\{D\})=dim\{D\}+1.$

For the purpose of convenience we introduce the concept of $L$-base of linear
system in the following way. Suppose $\{D\}$ is a linear system of Weil
divisors, $\{L\}$ is a family of curves parameterized by base $S.$ For every
nonempty Zariski open subset $U\subset S$ let $V(U,\{D\})$ be a linear subspace
in $H^0(\{D\})$, spanned by $s,$ such that $(s)$ contains $L_u$ for some $u\in
U.$ Evidently, $V(U^\prime\bigcap U^{\prime \prime},\{D\})\subset
V(U^\prime,\{D\})\bigcap V(U^{\prime \prime},\{D\})$ and $H^0(\{D\})$ is
finite-dimensional. Therefore there exists the minimal $V(U^*,\{D\}),$ such
that $V(U^*,\{D\})\subset V(U,\{D\})$ for every $U\subset S.$ Then
$|V(U^*,\{D\})|$ will be called $L$-base of $\{D\}$ and denoted by $\{D\}^L.$
\begin{Proposition}\label{S6}
$h^0(\{D\}^L)\ge h^0(\{D\})-L\cdot D-1$
\end{Proposition}
{\bf Proof} Suppose $\{D\}^L=|V(U^*,\{D\})|,$ $u\in U^*.$ We can also assume
that $L_u$ is not contained  in $Sing(X).$ Choose on $L_u$ points $x_1,$ $x_2,$
. . . , $x_d,$ $L\cdot D<d\le L\cdot D+1$ lying in nonsingular part of $X.$ The
condition of vanishing in  $x_1,$ $x_2,$ . . . , $x_d$ determines a subspace in
$H^0(\{D\})$ of codimension no greater than $d$ and $d\le L\cdot D+1.$ Now we
just notice that for every $s$ from this subspace $(s)$ contains $L_u$, because
otherwise we would have a contradiction by intersecting it with $L_u$.

 $ $

Define a new linear system $\{H_*\}$ by the following procedure. Denote $|H|$
by $\{H_0\}$ and for every nonnegative integer $i$ let $\{H_{i+1}\}$ be a
movable part of  $\{H_i\}^L.$ Evidently, $\{H_i\}$ will eventually stabilize.
This stabilized $\{H_i\}$ will be our  $\{H_*\}.$ It is evident that  $\{H_*\}$
is movable and  $\{H_*\}=\{H_*\}^L.$ (Here we set as definition that trivial
linear systems $\emptyset$ and $|O_X|$ are movable.)
\begin{Proposition}\label{S6*}
If $h^0(H)>1+(c_1+1)(c_2+c_1c_4+1)$ then $\{H_*\}$ is not trivial
\end{Proposition}
{\bf Proof} First of all, let $\{H\}^L=a_iE+D_i+\{H_{i+1}\}$, where $a_i\ge 0$,
$D_i$ does not contain E. Notice that if $a_i=0$ then
$\{H_{i+1}\}^L=\{H_{i+1}\}$ and the procedure stabilizes. From the other hand,
$\sum a_i\le c_1$ because $E\cdot L\ge 1$ and $H\cdot l\le c_1.$ Therefore
$\{H_*\}=\{H_{[c_1]+1}\}.$ It is easy to see that for all $i$ $H_i\cdot L\le
H\cdot L+c_1(-E\cdot L)\le c_2+c_1c_4.$ Therefore by proposition \ref{S6} we
have that $h^0(\{H_*\})\ge h^0(H)-(c_1+1)(c_2+c_1c_4+1)>1$. This implies
$\{H_*\}$ is not trivial.

 $ $

We also have from the above proof that $H_*\cdot L\le c_2+c_1c_4.$ Apply to
$K_X+2\{H_*\}$ Alexeev Minimal Model Program (\cite{Al}). Namely, let $\pi
:Y_1\longrightarrow X$ be a terminal modification of $K_X+2\{H_*\}$ in sense of
Alexeev.
\begin{Proposition}\label{S7}
Under the above notations the following is true.

$(1)$ $Y_1$ is Q-factorial and have at worst terminal singularities.

$(2)$ $\{\pi ^\prime H_*\}$ is free. Here $\{\pi ^\prime H_*\}$ is a inverse
image of linear system $\{H\}$ in sense of Alexeev, that is general element of
$\{\pi ^\prime H_*\}$ is $\pi ^\prime H_*$ for general $H_*\in \{H\}.$

$(3)$ $K_{Y_1}\cdot L\le c_3+2(c_2+c_1c_4)$
\end{Proposition}
{\bf Proof} Parts $(1)$ and $(2)$ are proved the same way as lemma $1.22$ in
\cite{Al}. Part $(3)$ is a corollary of the following chain of inequalities.

$K_{Y_1}\cdot L\le (K_{Y_1}+2(\pi ^\prime H_*))\cdot L \le (K_X+2H_*)\cdot L\le
c_3+2(c_2+c_1c_4)$

Here the middle inequality is due to the following argument. By definition of
terminal modification $K_{Y_1}+2(\pi ^\prime H_*)$ is $\pi -nef$ and therefore
in adjunction formula $K_{Y_1}+2(\pi ^\prime H_*)=\pi ^*(K_X+2H_*)+\sum
a_iD_i,$ where $D_i$ are exceptional divisors, all $a_i\le 0.$

 $ $

For the rest of the section we will use the following notations. Suppose $D_i,$
$i=1,...,k$ are exceptional divisors of morphism $\pi .$ For an arbitrary Weil
divisor $F$ on $X$ we will say that discrepancy of $F$ is a $k$-tuple
$\{discr_{D_i}(F)\}$ of discrepancies of $F$ in $D_i$, that is numbers
$discr_{D_i}(F)$ from the formula $\pi ^*F=\pi ^\prime (F)+\sum
discr_{D_i}(F)D_i.$
In these notations we have the following lemma.
\begin{Lemma}\label{L5}
Suppose $F=(s), s\in H^0(O_X(F))$. Suppose $s=\sum \alpha _j s_j,$ where
$(s_j)=F_j.$ Then for all $D_i$ $discr_{D_i}(F)\ge \min{_j}discr_{D_i}(F_j)$
and for a general $\{\alpha _j\}$ for given $\{s_j\}$ this inequality becomes
an equality.
\end{Lemma}
{\bf Proof} Suppose $rF$ is a Cartier divisor. In a neighborhood of generic
point $\pi (D_i)$ the sheaf $O_X(rF)$ can be trivialized. With respect to this
trivialization the local equation $f$ of divisor $rF$ is, by Newton binomial
formula, a linear combination of local equations $f_{(\gamma)}$ of divisors
$\sum \gamma _jF_j$, where $\sum \gamma _j=r,$  $\gamma _j\in {\bf Z}_{\ge 0}.$
By the definition, $discr_{D_i}(F)=\frac1r discr_{D_i}(rF)$ and
$discr_{D_i}(rF)$ is just an image of $f$ by a valuation on function field $\bf
C(X)$ of variety $X$ corresponding to $D_i.$ Therefore for arbitrary $\{\alpha
_j\}$  $discr_{D_i}(rF)\ge \frac1r\min{_j}discr_{D_i}(\sum \gamma _jF_j)\ge
\min{_j}discr_{D_i}(F_j)$ and for general $\{\alpha_j\}$ it becomes an
equality.

 $ $

Suppose now that $P_1\subset \{H_*\}$ is a set of all divisors $H_*$ containing
some $L\in \{L\}.$ Suppose a general element of $P_1$ has discrepancy
$\{d_i\}.$ Denote the set of all divisors from $P_1$ with such discrepancy by
$P.$
\begin{Proposition}\label{S9}
"Equations" of $H_*$, $H_*\in P,$ span $H^0(\{H_8\}).$
\end{Proposition}
{\bf Proof} For a general $L\in \{L\}$ divisors $H_*\in P,$ containing $L$
constitute a nonempty Zariski open subset in linear system of divisors from
$\{H\}$ containing $L.$ Therefore their "equations" span the corresponding
subspace in $H^0(\{H_*\}).$ By definition $\{H_*\}=\{H_*\}^L,$ so we are done.
\begin{Proposition}\label{S8}
$\{L\}$ have no base points on $E_{Y_1}.$
\end{Proposition}
{\bf Proof}
Proposition \ref{S9} and lemma \ref{L5} applied together imply that discrepancy
of general element of linear system $\{H_*\}$ equals $\{d_i\}.$ Therefore for
every $H_*\in P$ $\pi ^\prime H_*\in \{\pi ^\prime H_*\}.$ Moreover, the linear
equivalence between divisors $\pi ^\prime H_*$ is given by the same functions
from $\bf C(Y_1)=\bf C(X)$ as between corresponding divisors $H_*.$ Therefore
the proposition \ref{S9} implies that "equations" of  $\pi ^\prime H_*,$ where
$H_*\in P,$ span $H^0(\{\pi ^\prime H_*\}).$.

Suppose all $L$ on $Y_1$ pass through some point $y.$ Then all $\pi ^\prime
H_*,$ where $H_*\in P,$ contain $y$. But it is in  contradiction with
proposition \ref{S7}, $(2),$ so proposition \ref{S8} is proven.

 $ $

To complete the proof of the whole Accurate Resolution Lemma it is enough to
choose an arbitrary resolution of singularities $Y\longrightarrow Y_1.$ Then
$Y\longrightarrow X$ will satisfy all the requirements of accurate resolution.

\section{Treatment of case (2)}
Now we are in situation and notations of case $(2).$ (See section $3$.)
\begin{Proposition}\label{S10}
On $Y_1$ there exists a divisor $E$ which is exceptional with respect to
morphism $\pi^{Y_1}_X:Y_1\longrightarrow X$ such that $E\cdot l\ge 1.$
\end{Proposition}
{\bf Proof} Suppose $C$ is some general enough curve on the image $Z$ of
RC-fibration $\phi .$ Suppose $D\subset X$ is an image by $\pi^{Y_1}_X$ of the
surface $[ {\phi^{-1}(C)}].$ (Here parenthesis means Zariski closure.) The
general $l_{Y_1}$ does not intersect $\phi^{-1}(C)$ and, therefore,
$[{\phi^{-1}(C)}].$ ($Y_1$ is smooth therefore $\{l\}$ is free, see \cite{Na}.)
So, if $l_{Y_1}$ does not intersect with exceptional divisors of $\pi^{Y_1}_X$
then $l_X\cdot D=0,$ that is impossible because $X$ is Q-factorial and $\rho
(X)=1.$ Q.E.D.

 $ $

Notice that if $E\cdot l\ge 1$ then general $l_{Y_1}$ intersects with $E$ in
general points because $\{l_{Y_1}\}$ is free. Two cases are possible.

$(A)$ There exists such $E\subset Y_1$ that is exceptional with respect to the
morphism $\pi ^{Y_1}_Y.$

$(B)$ Family  $\{l_Y\}$ is free. Then there exists $E\subset Y$ that is
exceptional with respect to $\pi ^Y_X.$

The proof is generally the same in both cases but some technical details are
different. We begin with the case $(A).$ By the relative version of the usual
Minimal Model Program morphism $\pi ^{Y_1}_Y$ can be decomposed into extremal
contractions and flips, relative over $Y.$ Suppose $\pi ^{Y_3}_{Y_2}$ is the
first that contract some divisor $E_{Y_3},$ for which $l\cdot E_{Y_3}\ge 1.$
Suppose $\hat{E_{Y_3}}$ is a minimal resolution of $E_{Y_3}.$
\begin{Proposition}\label{11a} {\bf (Case (A))}
There exists a covering family $\{L\}$ of rational curves on $E_{Y_3},$ such
that the following conditions hold true.

$(1)$ $L\cdot K_{Y_3}<0$

$(2)$ $-L\cdot E_{Y_3}<3$

$(3)$ $L$ does not admit a nontrivial 2-point deformation on $\hat E_{Y_3}$,
that is a deformation with two fixed points, whose image is not in $L.$
\end{Proposition}
{\bf Proof} Suppose $\pi ^{Y_3}_{Y_2}(E_{Y_3})$ is a curve. Then we can choose
$\{L\}$ to be the fibers of $\pi ^{Y_3}_{Y_2} |_{E_{Y_3}}.$ Then $(1)$ is true
by the definition of extremal contraction. Suppose $\tilde E_{Y_3}$ is a
normalization of $E_{Y_3}.$ Then $\{L\}$ does not have base points on $\tilde
E_{Y_3}$ and therefore $L$ does not pass through its singularieties. This
easily implies $(3).$ The condition $(2)$ follows from the fact that (by lemma
\ref{LA}) $$(K_{Y_3}+E_{Y_3})\cdot L\ge K_{\hat E_{Y_3}}\cdot L=-2>-3.$$
Suppose now that $\pi ^{Y_3}_{Y_2}(E_{Y_3})$ is a point. Consider a minimal
model $F$ of $\hat E_{Y_3}.$ The surface $\hat E_{Y_3}$ is birationally ruled
or rational therefore we have two possibilities for $F$:

$1)$ $F\cong P^2$

$2)$ $F$ is ruled, there is a morphism $\theta :F\longrightarrow C$

We let $\{L\}$ be the family of planes on $P^2$ in the first case and the
family of fibers of $\theta$ in the second one. It evidently satisfies the
condition $(3).$ The condition $(1)$ holds for arbitrary curve on $E_{Y_3}$.
The condition $(2)$ again follows from the fact that $$(K_{Y_3}+E_{Y_3})\cdot
L\ge K_{\hat E_{Y_3}}\cdot L\ge -3.$$ The proposition is proven.

 $ $

Now we can apply the Accurate Resolution Lemma (lemma \ref{LR}.) Here $X$ means
$Y_3$, $H$ means $(\pi ^{Y_3}_X)^*(-2nK_X)$ and constants will be as follows.

$c_1=12n,$ $c_2=0,$ $c_3=0,$ $c_4=3.$

We see that if $h^0(-2nK_X)$ is big enough there exists a resolution
$Y_4\longrightarrow Y_3$ such that $K_{Y_4}\cdot L\le 2(3\cdot 12n)=72n$ and
$\{L\}$ have no base points on $E_{Y_4}.$
\begin{Proposition}\label{S12}
$L$ does not admit a nontrivial 2-point deformation on $Y_4.$
\end{Proposition}
{\bf Proof} If such deformation existed it would be a deformation on $E_{Y_4}$
by rigidity lemma. (About this lemma see \cite{CKM}, section 1. I must only
notice that it is not stated there correctly, one should add a condition of
flatness of morphism $f.$ It was noticed by several people, my attention was
brought to it by Iskovskikh.) The system $\{L\}$ has no base points on
$E_{Y_4}$ therefore $L$ does not pass through the singularities of
normalization $\tilde{E_{Y_4}}$ of the surface $E_{Y_4}.$ Resolution of
singularities $\hat{E_{Y_4}}$ is naturally mapped to $\hat {E_{Y_3}}$ therefore
2-point deformation of $L$ on $E_{Y_4}$ gives deformation on $\tilde{E_{Y_4}}$
and then on $\hat{E_{Y_4}},$ and then on $\hat {E_{Y_3}}.$ The last is
impossible by the choice of $L$, Q.E.D.

 $ $

Now we can apply to $\{L\}$ and $\{l\}$ the gluing lemma on $Y_4$ (see
\cite{KoMiMo1}) to obtain a new covering family of rational curves $\{l^\prime
\}.$ But now the image of RC-fibration corresponding to $\{l^\prime \}$ has
dimension $1$ or $0.$ And $l^\prime \cdot (-K_X) \le (1+dimY_4+L\cdot
K_{Y_4})(l\cdot (-K_X))\le 6(4+72n).$ So we managed to reduce the case $(2A)$
to cases $(1)$ and $(0),$ as it was promised at the end of section 3.

 $ $

Now we consider the case $(B).$ Similarly to the case $A,$ we have the
following statement.
\begin{Proposition}\label{11b} {\bf (Case (B))}
There exists a covering family $\{L\}$ of rational curves on $E_{Y},$ such that
the following conditions hold true.

$(1)$ $L\cdot K_{Y}<0$

$(2)$ $-L\cdot E_{Y}<3$

$(3)$ $L$ does not admit a 2-point nontrivial deformation on $\hat E_{Y_3}.$

$(4)$ $\pi ^Y_X(L)$ is a point
\end{Proposition}
{\bf Proof} If $\pi ^Y_X(E_Y)$ is a curve let $\{L\}$ be the family of fibers
of $\pi ^Y_X|_{E_Y}.$ If $\pi ^Y_X(E_Y)$ is a point then let it come from the
minimal model of $\hat E_Y$ as in the proof of proposition \ref{11a}. As in the
case $(A),$ $K_{\hat E_Y}\cdot L$ is $-2$ or $-3.$ Conditions $(3)$ and $(4)$
are evidently satisfied, we only need to prove $(1)$ and $(2).$ In order to do
it consider the adjunction formula for $\pi ^Y_X,$ multiplied by $L:$
$$K_Y\cdot L=\sum_{E_i\ne E_Y}a_iE_iL+aE_Y\cdot L
{}~~~~~~~~~~~~~~~~~~~~~~~~~~~~~~(*) $$
Here $a_i$ and $a$ are discrepancies, they are of form $(-\frac mn),$ $m\in
\{0,1,...,n-1\},$ where $n$ is an index of $X.$ (Discrepancies are nonpositive
because $Y$ is a terminal modification of $X.$) We have the following chain of
inequalities.
$$-3\le K_{\hat E_Y}\cdot L\le (1+a)E_Y\cdot L+\sum_{E_i\ne E_Y}a_iE_iL\le
(1+a)E_Y\cdot L$$
Here the middle inequality follows from lemma \ref{LA} and formula $(*),$ and
the right from nonpositivity of $a_i.$ Therefore $1+a\ge \frac1n$ implies that
either $-E_Y\cdot L\le 0$ or $-E_Y\cdot L\le 3n.$ Therefore $-E_Y\cdot L\le
3n.$
Now the condition $(1)$ follows from the following chain of inequalities.
$$K_Y\cdot L=\sum_{E_i\ne E_Y}a_iE_iL+aE_Y\cdot L\le aE_YL\le 3n$$
Here the right inequality holds because of the following argument. We know that
$-1<a\le 0$ therefore $E_YL\ge 0$ implies $aE_YL\le 0$ and  $E_YL< 0$ implies
$aE_YL\le -E_YL.$ Q.E.D.

 $ $

Again, as in case $(A),$ we apply the Accurate Resolution Lemma (lemma
\ref{LR}.) The only difference is that now we have $Y$ instead of $Y_3$ and
constants are as follows.

$c_1=12n,$ $c_2=0,$ $c_3=3n,$ $c4=3n.$

Again if $h^O(-2nK_X)$ is big enough there exists an accurate resolution $Y_4.$
 We have again that $L$ does not admit nontrivial 2-point deformation on $Y_4$.
(Arguments from the proof of proposition \ref{S12} work without any problems
because of condition $(4)$ of proposition \ref{11b}.) So we can apply gluing
lemma from \cite{KoMiMo1}. The bound on $l^\prime \cdot (-K_X)$ will be the
following.

$l^\prime \cdot (-K_X) \le (4+L\cdot K_{Y_4})(l\cdot (-K_X))\le
(4+3n+2(12n\cdot 3n))\cdot 12n=12n(4+3n+72n^2).$

So we completed the treatment of case $(2B)$. Our Main Theorem is finally
proven.

\end{document}